\documentclass{article}
\usepackage{amssymb,amsmath,amsthm}
\newtheorem{theorem}{Theorem}

\newtheorem{lemma}{Lemma}

\date{}
\numberwithin{equation}{section} \numberwithin{theorem}{section}
\numberwithin{lemma}{section} \numberwithin{corollary}{section}
\numberwithin{remark}{section} \numberwithin{proposition}{section}
\numberwithin{definition}{section}
\begin{document}
% let's define a new macro
\newcommand{\n}{\noindent}
\newcommand{\vs}{\vskip}
\title{A Note on Lipschitz Continuiy of the Solutions 
of a Class of Elliptic Free Boundary Problems}

\vskip 0.5cm
\author{A. Lyaghfouri \\
American University of Ras Al Khaimah\\
Department of Mathematics and Natural Sciences\\
Ras Al Khaimah, UAE}
\maketitle

 \vskip 0.5cm
\begin{abstract}
We provide a new and simple proof based on Harnack's inequality to the Lipschitz
continuity of the solutions of a class of free boundary problems.
\end{abstract}
\vskip 0.5cm

\section{Introduction}\label{1}

\vskip 0.5cm

We consider the following problem
\begin{equation*}(P)
\begin{cases}
& \text{ Find }\,\,(u, \chi) \in  H^1(\Omega)\times L^\infty
(\Omega)
\text{ such that}:\\
&(i)\quad u=0\,\,\text{ on } \Gamma\\
& (ii)\quad  0\leqslant u\leqslant M, \quad 0\leqslant\chi\leqslant
1 ,
\quad u(1-\chi) = 0 \,\,\text{ a.e.  in } \Omega\\
& (iii)\quad div(a(x)\nabla u+\chi H(x))=0 \quad\text{in }
 \mathcal{D}'(\Omega),
\end{cases}
\end{equation*}
where $\Omega$ is an open bounded domain of $\mathbb{R}^n$,
$\Gamma$ is an open nonempty $C^{1,1}$ subset of
the boundary of $\Omega$, $M$ is a positive number, $a(x)= (a_{ij}(x))$ is an $n$-by-$n$ matrix
that satisfies for some positive constants $\lambda$, $\Lambda$, and $\alpha\in(0,1)$
\begin{eqnarray}\label{e1.1-1.5}
&a\in C^{0,\alpha}_{loc}(\Omega\cup\Gamma),\\
&\sum_{i,j}|a_{ij}(x)|\leq \Lambda,\quad\text{  for all }x\in \Omega,\\
&{ a}(x)\xi.\xi \geq \lambda\vert \xi\vert^2,\quad\text{ for all }(x,\xi)\in
\Omega\times\mathbb{R}^{2}.
\end{eqnarray}

\n $H(x)$ is a vector function satisfying for positive constants
$\bar h $ and $p>n/(1-\alpha)$
\begin{eqnarray}\label{2.1-2}
 & |H|_\infty\leqslant \overline{h},\\
&div(H) \in L_{loc}^p(\Omega\cup\Gamma).
\end{eqnarray}

\n Problem $(P)$ describes various free boundary problems including the
heterogeneous dam problem \cite{[A1]} \cite{[ChiL1]}, \cite{[CL1]}, \cite{[CL3]}, 
\cite{[L1]}, \cite{[L2]}, \cite{[L3]},
in which case $\Omega$ represents a porous medium with permeability matrix $a(x)$,
and $H(x)=a(x)e$ with $e=(0,...,0,1)$.
The weak formulation of the lubrication problem \cite{[AC]} is
obtained for $a(x)=h^3(x)I_2$ and $H(x)=h(x)e$, where $I_2$ is the
2-by-2 identity matrix,
and $h(x)$ is a scalar function related to the Reynolds equation.
A third problem is the aluminium electrolysis \cite{[BMQ]} obtained 
for $a(x)=k(x)I_2$ and $H(x)=h(x)e$, where $k(x)$ and $h(x)$ are scalar 
functions.

\vs 0.2cm\n In this problem the free boundary $\partial\{u>0\}\cap\Omega$ represents the
interface between the two sets $\{u=0\}$ and $\overline{\{u>0\}}$.
For the dam and lubrication problems, the free boundary is the
interface that separates the region containing the fluid from the rest of the domain.
For the aluminium electrolysis problem, it represents the region containing liquid
and solid aluminium.

\vs 0.5cm\n Since $\chi H(x)\in L_{loc}^\infty(\Omega)$ and due
to (1.1)-(1.5), it follows from $(P)iii)$, \cite{[GT]} Theorem 8.24, p. 202
and Theorem 8.29, p. 205 that $u\in C_{loc}^{0,\beta}( \Omega)$ for some $\beta\in(0,1)$.
In this paper we will show that $u\in C_{loc}^{0,1}( \Omega)$.
We observe that Lipschitz continuity is the optimal regularity because of
the gradient jump over the free boundary. Lipschitz continuity is important
due to its major role in the free boundary regularity as
in \cite{[CL2]} and \cite{[CL4]} for example.

\section{Interior Lipschitz Continuity}\label{2}

\vs 0,5cm \n Under the above assumptions, we have the following
interior regularity result.

\begin{theorem}\label{t2.1}
Let $(u,\chi)$ be a solution of the problem $(P)$. 
Then we have $~u\in C^{0,1}_{loc}(\Omega).$
\end{theorem}

\vs0.3cm \n Theorem 2.1 was previously established in \cite{[CL2]} and \cite{[L4]}.
The same method was successfully extended to the quasilinear case in
\cite{[CL5]} and \cite{[CL6]}.
Here we would like to propose a different approach based on the Harnack inequality
which does neither require $div(a(x)(x-y))$ to be
uniformly bounded in $y$ from above nor that $div(H)$ be uniformly bounded 
from below as in \cite{[CL2]} and \cite{[L4]}.

\vs 0.3cm\n We need a lemma.

\begin{lemma}\label{l2.2} Let $x_0\in \Omega$ and $r>0$ such that
$B_r(x_0)\subset \{u>0\}$, $B_{5r}(x_0)\subset\subset \Omega$,
and $\partial B_r(x_0)\cap \partial \{u>0\}\neq \emptyset$. Then we
have for some positive constant $C$ depending only on $n$, $p$,
$\lambda$, $\Lambda$ and $\bar h$
$$\displaystyle{\max_{\overline{ B_{r}}(x_0)} u \,\leq \,C
\,r}.$$
\end{lemma}

\vs 0,3cm\n \emph{Proof.} Let $x_1\in \partial B_r(x_0)\cap \partial \{u>0\}$.
First since $B_{5r}(x_0)\subset\subset \Omega$, it is easy to verify that
$B_{2r}(x_1)\subset\subset \Omega$.
Next we apply the Harnack inequality \cite{[GT]} p. 194 to the equation
$(P)iii)$, and we get for a positive constant $C$ depending only on $n$, 
${\Lambda\over\lambda}$, and $p$
\begin{eqnarray*}
\max_{\overline{B}_{2r}(x_1)} ~ u&\leq& C
\Big(\min_{\overline{B}_{2r}(x_1)} ~ u +
{1\over\lambda}r^{1-{n\over p}}|\chi H|_{p,\overline{B}_{2r}(x_1)}\Big)\nonumber\\
&\leq& C\Big(0+{\overline{h}\over\lambda}(2r)^{1-{n\over p}}.|B_{2r}(x_1)|^{1\over p}\Big)\nonumber\\
&=& \frac{C\overline{h}|B_1|^{1\over p}}{\lambda}.(2r)^{1-{n\over p}}.(2r)^{{n\over p}}\nonumber\\
&=& \frac{2C\overline{h}|B_1|^{1\over p}}{\lambda}r.
\end{eqnarray*}

\n Given that $B_r(x_0)\subset B_{2r}(x_1)$, the lemma follows.
\qed

\vs 0.3 cm\n \emph{Proof of Theorem 2.1}. Using Lemma 2.1,
the proof follows as in \cite{[L4]}.
\qed

\section{Boundary Lipschitz Continuity}\label{3}

\vs 0,5cm \n Under the assumptions (1.1)-(1.5), we have the following
boundary regularity result.

\begin{theorem}\label{t3.1}
Let $(u,\chi)$ be a solution of the problem $(P)$. 
Then we have $~u\in C^{0,1}_{loc}(\Omega\cup\Gamma).$
\end{theorem}

\vs 0.3cm\n We need a lemma.

\begin{lemma}\label{l3.1}
It is enough to establish Theorem 3.1 when $\Gamma$ is part of a hyperplane.
\end{lemma}

\vs 0,3cm\n \emph{Proof.} Indeed let $x_0\in \Gamma$.
Since $\Gamma$ is a $C^{1,1}$ manifold, there exists
an open set $U\subset \mathbb{R}^n$ that contains $x_0$
and a $C^{1,1}-$diffeomorphism $\Upsilon:~B_1~\rightarrow~U$ such that
\begin{equation}\label{3.1}
\Upsilon(B_1\cap\{y_n<0\})=U\cap\Omega\quad\text{and}\quad \Upsilon(B_1\cap\{y_n=0\})=U\cap\partial\Omega,
\end{equation}
where $B_1$ is the unit ball of $\mathbb{R}^n$.

\n Let $J\Upsilon(y)$ be the Jacobian determinant of $\Upsilon$, and set
$\Sigma=B_1\cap\{y_n=0\}$, $B_1^-=B_1\cap\{y_n<0\}$, $v(y)=uo\Upsilon(y)$, $\theta(y)=\chi o\Upsilon(y)$,
$b(y)=|J\Upsilon(y)|.^t(D\Upsilon^{-1}(\Upsilon(y)))a(\Upsilon(y))D\Upsilon^{-1}(\Upsilon(y))$,
and $K(y)=|J\Upsilon(y)|.^t(D\Upsilon^{-1}(\Upsilon(y))) Ko\Upsilon(y)$. Then 
it is not difficult to verify that $(v,\theta)$ is a solution of the following problem:

\begin{equation*}(Q)
\begin{cases}
& \text{ Find }\,\,(v, \theta) \in  H^{1}(B_1^-)\times L^\infty (B_1^-)
\text{ such that}:\\
&(i)\quad v=0\,\,\text{ on } \Sigma\\
& (ii)\quad  0\leqslant v\leqslant M, \quad 0\leqslant \theta\leqslant
1 ,
\quad v(1-\theta) = 0 \,\,\text{ a.e.  in } B_1^-\\
& (iii)\quad div(b(y)\nabla v+\theta K(y))=0 \quad\text{in }
 \mathcal{D}'(B_1^-)
\end{cases}
\end{equation*}

\n Moreover, it is clear that $b$ and $K$ satisfy the assumptions (1.1)-(1.5) in $B_1^-$.
\qed

\vs 0.3cm\n Henceforth we will establish Theorem 3.1 for the solutions of the problem $(Q)$.
To do that we need a lemma.

\begin{lemma}\label{l3.2}
For each $y_0\in \sum$ and each ball $B_{2r}(y_0)\subset\subset B_1$, there
exists a positive constant $C$ such that
$$\max_{\overline{B}_r^-(y_0)} ~ v \leq Cr.$$
\end{lemma}

\vs 0.3 cm\n \emph{Proof}. We will use the reflection method and introduce the following extensions
to $B_1$ of a function $f$ defined in $B_1^-$
\begin{equation*}
\widetilde{f}(y)=\left\{
\begin{array}{ll}
f(y),  &\hbox{ if } y\in B_1^-\\
f(y',-y_n),  &\hbox{ if } y\in B_1\setminus B_1^-
       \end{array}
     \right.
~\text{ and }~
f^*(y)=\left\{
\begin{array}{ll}
f(y),  &\hbox{ if } y\in B_1^-\\
-f(y',-y_n), &\hbox{ if } y\in B_1\setminus B_1^-.
\end{array}
\right.
\end{equation*}

\n Then we have by arguing as in \cite{[CL3]}
\begin{equation*}
\left\{
  \begin{array}{ll}
& div(c(y)\nabla \widetilde{v} +\widetilde{\theta} G(y))=0\qquad\text{in}\quad
{\cal D}'(B_1)\\
& \widetilde{v}\geq 0,\quad 0\leq\widetilde{\theta}\leq 1,\quad \widetilde{v}(1-\widetilde{\theta})=0, \qquad\text{in}\quad B_1.
  \end{array}
\right.
\end{equation*}

\n where $c(y)=(c_{ij}(y))$ is the $n$-by-$n$ matrix defined by

\begin{equation*}
c_{ij}=\left\{
\begin{array}{ll}
\widetilde{b}_{ij}, & \hbox{if } 1\leq i,j\leq n-1  \hbox{ or } i=j= n\\
b^*_{ij}, & \hbox{if } 1\leq i\leq n-1\hbox{ and }j=n \hbox{ or } i=n \hbox{ and } 1\leq j\leq n-1.
       \end{array}
     \right.
\end{equation*}

\n and $G(y)=(G_1(y),...,G_n(y))$ is the vector function defined by
\begin{equation*}
G_i=\left\{
\begin{array}{ll}K^*_i, & \hbox{if } 1\leq i\leq n-1\\
\widetilde{K}_n, & \hbox{if } i=n.
       \end{array}
     \right.
\end{equation*}

\n We observe that the matrix $c(y)$ satisfies (1.2)-(1.3) and that $G$ satisfies (1.4)
in $B_1$. As in the proof of Lemma 2.1, we get by Harnack's inequality, for some positive 
constant $C$.
\begin{eqnarray}\label{e3.2}
\max_{\overline{B}_r(y_0)} ~ \widetilde{v}&\leq& C r.
\end{eqnarray}
\n The lemma follows from (3.2).
\qed

\vs 0.3 cm\n \emph{Proof of Theorem 3.1}. Using Lemmas 3.1 and 3.2,
the proof follows as in \cite{[L4]}.
\qed

\end{document}